\RequirePackage{fix-cm}
\documentclass[smallextended]{svjour3}       
\smartqed  
\usepackage{lineno,hyperref}
\modulolinenumbers[5]
\usepackage{amssymb,latexsym}
\usepackage{enumerate}
\usepackage{amsmath}
\usepackage{graphicx}
\usepackage{verbatim}
\usepackage[all]{xy}
\usepackage{epsfig}
\usepackage{mathrsfs}
\usepackage[utf8]{inputenc}
\usepackage{calc} 
\usepackage{hypcap}
\usepackage[font=footnotesize,labelfont=bf]{caption}
\usepackage[font=scriptsize,labelfont=bf]{subcaption}
\usepackage[mathscr]{eucal}
\usepackage{bbold}
\usepackage{dsfont}
\usepackage{graphics}
\usepackage{framed}
\usepackage{fancyhdr}
\usepackage{footnote}
\usepackage{bbm}
\usepackage{float,lscape}

\makesavenoteenv{equation}
\newcommand{\Pds}{\mathds{P}}
\newcommand{\Esf}{\mathsf{E}}

\newcommand{\thetahat}{\hat{\theta}}

\newcommand{\phitil}{\widetilde{\phi}}
\newcommand{\vhat}{\widehat{v}}
\newcommand{\Xtil}{\widetilde{X}}
\newcommand{\Ttil}{\widetilde{T}}
\newcommand{\Lcal}{\mathcal{L}}
\newcommand{\Khat}{\widehat{K}}

\newcommand{\Lcalhat}{\widehat{\mathcal{L}}}

\newcommand{\set}[1]{\left\{#1\right\}} 
\newcommand{\abs}[1]{\left|#1\right|}
\newcommand{\norm}[1]{\left\lVert #1\right \rVert}

\begin{document}

\title{Distribution free goodness of fit tests for regularly varying tail distributions 
}


\author{Thuong\ T.\ M.\ Nguyen}


\institute{Thuong T. M. Nguyen \at
              School of Mathematics and Statistics,
              Victoria University of Wellington\\ Wellington 6140, New Zealand \\
              \email{thuong.nguyen@vuw.ac.nz}           
           }


\maketitle

\begin{abstract}
We discuss in this paper a possibility of constructing a whole class of asymptotic distribution-free tests for testing regularly varying tail distributions. The idea is that we treat the tails of distributions as members of a parametric family and using MLE to estimate the exponent. No matter what the exponent's estimator is, we are able to transform the whole class into a specific distribution with a prefix exponent so that we are free from choosing any functional of the tail empirical process as a distribution-free test statistic. The asymptotic behaviour of some new tests, as examples from the whole class of new tests, are demonstrated as well.

\end{abstract}

\section{Introduction}
Suppose that we have a random sample $X_1,X_2,\cdots, X_n$ of a random variable $X$ following some unknown distribution $F$. No matter how $F$ 
behaves, we are 
mainly interested in the right tail behavior of $F$, i.e., in $F(x)$ when $F(x) \to 1$. Among the class of heavy-tail distributions, regularly varying 
tail distributions has been attracted most attention due to its various applications. Briefly, a distribution $F$ is said to be {\it regularly varying 
in tail} 
if 
\begin{align}
\lim_{x\to\infty} \frac{1-F(tx)}{1-F(x)} = t^{-\theta}
\end{align}
for all $ t >0$ where $\theta >0$ is called the {\it index} or {\it exponent} of regular variation. By this 
definition, regularly varying tail distributions are often indicated as Pareto-type distributions. 

Regular variation of the tail of a distribution often appears as a natural condition in various theoretical results of probability theory. 
The most typical example is that it is the condition for the distribution of the partial maxima 
to belong to the domain of attraction of extreme value distributions. This could be found in various references, we refer to books of de Haan and Ferreira 
\cite{deHaanbook} and Resnick \cite{Resnick} among others.
Moreover, regular variation naturally arises as a common phenomenon in numerous practical case studies such as finance, insurance, physics, geology, 
hydrology and 
engineering, etc. Therefore, detecting and testing regular variation plays an important role in probabilistic areas as well as in practical applications. 

\medskip

The problem of estimating the exponent $\theta$ has a very rich literature. Among various studies, the Hill estimator is the most 
common 
one, see Hill 
\cite{Hill}. Recall that the Hill estimator of the exponent $\theta$ essentially takes on the following form
\begin{align}
 \hat{\theta}^H = \Big(\frac{1}{k} \sum_{i=1}^k \ln X_{i,n} - \ln X_{k,n}\Big)^{-1}
\end{align}
where $k = k(n)$ is a certain sample fraction, $k(n) \to \infty$ in an appropriate way and $X_{1,n} \le \cdots \le X_{n,n}$ denotes the upper order 
statistics 
of the sample $X_1,\cdots, X_n$. Many other estimators was proposed, most of them is also based on the upper order statistics and is not too difficult to 
compute. These include, to name a few, popular estimators introduced by Dekkers et. al. \cite{Dekkersetal}, de Haan and Resnick \cite{DeHaan}, 
Pickands \cite{Pickands}, Teugels \cite{Teugels} among others. The rate of the convergence of an estimator was also discussed, for example, in Hall and
Welsh \cite{Hall84}.

\medskip

In contrast to the numerous number of approaches to estimate the exponent $\theta$, using goodness of fit (gof) for 
testing regular variation has been addressed in a modest number of studies.  Not long ago, Beirlant et al. \cite{Beirlantetall} modified 
the Jackson statistic - which was originally proposed as a gof for testing exponentiality - for testing Pareto-type data. Koning and Peng \cite{Koning} 
examined the Kolmogorov-Smirnov, Berk-Jones and the estimated score tests and compare them in terms of Bahadur efficiency. In that paper, the Berk-Jones 
and estimated score tests, which are not based on the empirical process, were shown to perform better than the Kolmogorov-Smirnov test.

\medskip

Note that regularly varying tail distributions are members of the family of the generalized Pareto distribution (GPD), which is 
of the form
\begin{align}\label{gpd}
 F(x;k,\sigma) =\begin{cases} 1 - (1-\frac{kx}{\sigma})^{1/k}, \quad k\ne 0, \sigma >0\\
 1 - e^{-x/\sigma}, \qquad\quad\,\, k=0, \sigma >0  
 \end{cases}
\end{align}
where $k$ and $\sigma$ are shape and scale parameters and $x >0$. To test the fit of data to a GDP, there has been several studies such as Davison and 
Smith \cite{Davison}, Choulakian and Stephens \cite{Choulakian}. In these papers, the critical values for Cramer-von Mises 
statistic 
and 
Anderson-Darling statistic, which is a weighted Cramer-von Mises statistic, were given.

\medskip

In this paper, we will introduce a wide class of asymptotically
distribution free gof tests for testing regularly varying tail distributions. Our approach follows a 
new method introduced in Khmaladze \cite{Estate2}. In that paper, the author proposed a unitary transformation which enables us to create a class of 
gof tests for both simple and parametric hypothesis testing problems. The method for the latter will be adopted for our problem, which we will present in 
Section \ref{method}. Briefly speaking, we will use exactly the same transformation in \cite{Estate2} to derive a modification of the empirical process. 
Only this time, the assumption that $F$ belongs to a parametric family of distributions is no longer available and the fact is that the right tail just 
partially describes the distribution $F$. Hence, we can only consider the tail empirical process and transform it. The transformed tail empirical process, 
under the hypothesis of interest, possesses a limit in distribution asymptotically free from any underlying distribution $F$ as well as the unknown 
index/exponent of the regular variation. Therefore, any appropriate functionals of the transformed process could be used as asymptotically distribution 
free test statistics. Section \ref{review} will be devoted for the main literature of the empirical process based on the tail.

\medskip

Some simulation results will be documented in Section \ref{simulation}. Namely, we will take the Kolmogorov-Smirnov (KS), Cramer-von Mises ($\Omega^2$) 
and Anderson-Darling ($A^2$) 
tests as 
examples from the new class of asymptotically distribution free gof tests for demonstration. The asymptotically distribution free property of these test 
statistics will be illustrated for different choices of the original distribution $F$.

\section{Main results}\label{method}

\subsection{Overall review}\label{review}

Suppose that from a random sample $X_1, \cdots, X_n$ we are only considering excesses over a certain threshold $x_0$ which 
is considerably large. Let us denote the subsample of all observations 
exceeding $x_0$ by $\Xtil_1,\cdots,\Xtil_m$. Generally, the choice of $x_0$ as well as the sample fraction $m/n$ may require some educated guessing. Hill 
\cite{Hill} suggested $m$ be chosen as an adaptive, data-analytic basis.
Several methods for choosing $m$ or $x_0$ based on survey data could be found in Drees and Kaufmann \cite{Drees}, Danielsson et al. \cite{Danielsson} 
and 
Guillou and Hall \cite{Guillou}. Nevertheless, there has existed various research assuming that $m$ is known 
such that $m\to \infty$ as $n\to \infty$ but $m=o(n)$. With this assumption and some prior knowledge about the underlying distribution function 
$F$, Haeusler and Teugels \cite{Haeusler} derived a general condition which can be used to determine the optimal value $m$ explicitly. Hall \cite{Hall82} 
also considered $m$ having a deterministic value and a quite common estimate of the exponent was introduced. Throughout this 
paper, our approach will be based on the same such assumption, that also means $x_0$ fixed.

\medskip

Our main aim is to create a class of gof test for testing the hypothesis $H_0:$ ``$F$ is a regularly varying tail distribution" against the alternative 
$H_1:$ ``$F$ is not a regularly varying tail distribution". It is sensible that the exponent $\theta$ needs to be estimated based on the right tail only. 
Specifically, we only consider observed values $\set{\Xtil_1,\cdots,\Xtil_m}$, which now should be looked at as a sample of a different random 
variable, let say, $\Xtil$.

\medskip

Denote $\Ttil = \frac{\Xtil}{x_0}$ and $T =\Ttil -1$. Under the hypothesis of interest and by the definition of the regular variation, the survival 
distribution of $\Xtil$ conditional on the specified large value $x_0$ is
\begin{align}
\frac{\Pds\set{\Xtil\ge x_0t}}{\Pds\set{\Xtil\ge x_0}} = \frac{1-F(x_0t)}{1-F(x_0)} = t^{-\theta} + o(1), \quad t\ge 1,\, \theta\ge 0.
\end{align}
Then, the 
distribution of the positive continuous random variable $T$ under the null hypothesis is
\begin{align}
 H_{\theta}(t) = 1 - (1+t)^{-\theta}, \qquad t \ge 0.
\end{align}
Clearly, the density function is $h(t) = \theta(1+t)^{-(\theta+1)}.$ Viewing $T$ as a new random variable, we can define the tail empirical process in a 
similar way as of the standard empirical process, i.e., we have
\begin{align}
H_m(t) = \frac{1}{m} \sum_{i=1}^m 1_{\set{T_i \le t}}.
\end{align}

\noindent The tail parametric
empirical process $\vhat_{mH}$ is 
\begin{align}\label{empro}
 \vhat_{mH}(t) = \sqrt{m}[H_m(t) - H_{\thetahat_m}(t)].
\end{align}
where $\thetahat_m$ is an estimator of $\theta$ calculated from the subsample $\Xtil_1,\cdots,\Xtil_m$ (or $T_1,\cdots, T_m$). Assume that the true 
unknown exponent under the hypothesis $H_0$ is $\theta_0$. Let $\thetahat_m$ be the maximum likelihood estimator (MLE), then it is 
the solution of the equation
\begin{align}
\sum_{i=1}^m \frac{\partial \log h(T_i, \theta)}{\partial \theta} =0.
\end{align}
That yields
\begin{align}
\thetahat_m = \frac{m}{\sum_{i=1}^m \log(T_i +1)},
\end{align}
which actually coincides with the Hill estimator $\thetahat^H$. This estimator was proved to be consistent in the sense that 
\begin{align}
 \thetahat^H = \thetahat_m \stackrel{\Pds}{\to} \theta_0
\end{align}
under the condition that $m\to \infty$ as $n\to\infty$ such that $m/n \to 0$. The proof of this 
convergence 
could be found in Mason \cite{Mason}. Regarding the asymptotic normality of the estimator, we refer to Haeusler and Teugels \cite{Haeusler}
, Geluk et al. \cite{Geluketal.}, de Haan and Resnick \cite{deHaanRes} among various others.

We will spend few more lines here to review the property of the limit in distribution of the process $\vhat_{mH}$ because it is 
essential for the method we present below. For all the concepts and terminologies, we follow and keep the same as in \cite{Estate2}. For the sake of 
lucidity, we extract and represent the method only for our particular problem. 

Consider the space $\Lcal_2(H)$. Recall that a function $\phi$ is integrable with respect to $H$ if $\int_0^{\infty} \phi(s) H(ds) <\infty$ and square 
integrable if
$\int_0^{\infty} \phi^2(s) H(ds) <\infty$. The space $\Lcal_2(H)$ consists of all square integrable functions with respect to $H$. The inner product and 
norm in $\Lcal_2(H)$ are defined as usual. That is, for any $\phi,\phitil \in 
\Lcal_2(H)$  we have
\begin{align*}
 \norm{\phi}_H^2 &= \int_0^{\infty} \phi^2(s) H(ds),\\
 \langle \phi, \tilde{\phi}\rangle_H &= \int_{0}^{\infty} \phi(s) \tilde{\phi}(s)H(ds).
\end{align*}

Denote by $w_H(\phi)$ a function-parametric $H$-Brownian motion where $\phi$ is a square integrable function in 
$\Lcal_2(H)$. That 
means, 
$w_H(\phi)$ for each $\phi$ is a Gaussian random variable with expected value $0$ and variance $\Esf w^2_H(\phi) = 
\norm{\phi}_H^2.$ This also implies that the covariance between $w_H(\phi)$ and $w_H(\tilde{\phi})$ is

\begin{align*}
 \Esf w_H(\phi)w_H(\tilde{\phi}) = \langle \phi, \tilde{\phi}\rangle_H.
\end{align*}
If $\phi_t(t^*) = 1_{\set{t^*\le t}}$ then $w_H(\phi_t):= w_H(t)$ is simply the Brownian motion in time $H(t).$ As usual, a linear transformation of 
$w_H(t)$ which is of the form $v_H(t)=w_H(t) - H(t) w_H(\infty)$ is the Brownian bridge in time $H(t).$

Assuming that under the null hypothesis $H_0,$ the true unknown exponent of the regular variation is $\theta_0$. Denote by $v_{mH}$ the usual tail 
empirical process, that is,
\begin{align}
 v_{mH}(t) = \sqrt{m}[H_m(t) - H_{\theta_0}(t)].
\end{align}

\noindent Let us stress that, since we assumed $x_0$ is known, the limit in distribution of 
$v_{mH}(t)$, as a {\it conditional} tail empirical process, is 
similar to that of the standard empirical process, which is the Brownian bridge $v_H(t)$. For a full description on the property of general tail empirical 
process, we refer to Einmahl \cite{Einmahl90,Einmahl92}. It was proved in these papers that if such process is {\it unconditional} on $x_0$, the limit in 
distribution of 
$v_{mH}(t)$ is a Brownian motion in time 
$H(t)$. On a recent approach on distribution free gof test for testing tail copula by Can et al. \cite{Canetal}, the tail empirical process - constructed 
on tail copula, was mapped to a standard Brownian motion. Their construction 
based on the innovative martingale method in Khmaladze \cite{Estate1} which is known as the Khmaladze transformation. Our approach here use a 
different Khmaladze transformation in \cite{Estate2}, and so let us call it Khmaladze-2.

\medskip

Considers the function-parametric version of the tail empirical process
\begin{align}\label{func.em.pr}
v_{mH}(\phi) = \int_0^{\infty} \phi(t) v_{mH}(dt) = \frac{1}{\sqrt{m}} \sum_{i=1}^m [\phi(T_i) - \Esf \phi(T_i)]
\end{align}
where $\phi$ is a function in $\Lcal_2(H)$. In a similar way to achieve the limit of $v_{mH}(t)$, the limit in distribution of $v_{mH}(\phi)$ with some 
proper restriction on functions $\phi$ is 
called a 
function-parametric $H$-Brownian bridge. Specifically, 
\begin{align}\label{limitvH}
 v_H(\phi) = w_H(\phi) - \langle \phi, \mathbb{1} \rangle_H w_H(\mathbb{1}),
\end{align}
where $\mathbb{1}$ stands for the function identically equals to $1$.



\medskip

We can also say more about the asymptotic behavior of the tail parametric empirical process where $\thetahat_m$ is 
the MLE which is also the Hill's estimator. It is well-known since long ago that under some usual and mild constraints, MLE 
$\thetahat_m$ possesses the asymptotic property
\begin{align}\label{mle-es}
 \sqrt{m}(\thetahat_m -\theta_0) = \Gamma_{H}^{-1} \int_0^{\infty} \frac{\dot{h}_{\theta_0}(t)}{h_{\theta_0}(t)} v_{mH}(dt,\theta_0) + o_P(1), 
\qquad m \to \infty
\end{align}
where $h_{\theta_0}$ denotes the hypothetical density and $\dot{h}_{\theta_0}$ its derivatives in $\theta$. Denote by
\begin{align}
 \Gamma_{H} = \int_0^{\infty} \frac{\dot{h}_{\theta_0}^2(t)}{h_{\theta_0}^2(t)} H(dt)
\end{align}
the Fisher information. Then it is easy to check that $\Gamma_{H} = 1/\theta_0^2.$ As a consequence of \eqref{mle-es}, we can expand the process 
$\vhat_{mH}(t)$ as
\begin{align}\label{em.pr-expand}
 \vhat_{mH}(t) &= v_{mH}(t) - \sqrt{m}[H_{\thetahat_m}(t)-H_{\theta_0}(t)]\nonumber\\
 &= v_{mH}(t) - \int_0^t \frac{\dot{h}_{\theta_0}(t)}{h_{\theta_0}(t)} H_{\theta_0}(dt) \Gamma_{H}^{-1}  \int_0^{\infty} 
\frac{\dot{h}_{\theta_0}(t)}{h_{\theta_0}(t)} v_{mH}(dt,\theta_0) + o_P(1)\nonumber\\
&= v_{mH}(t) - \int_0^t \beta_H^T(t) H_{\theta_0}(dt) \int_0^{\infty} \beta_H(t) v_{mH}(dt,\theta_0)+ o_P(1)
\end{align}
where 
\begin{align}
 \beta_H(t) = \Gamma_{H}^{-1/2} \frac{\dot{h}_{\theta_0}(t)}{h_{\theta_0}(t)}
\end{align}
denotes the normalized score function. This expression represents the limit in distribution $\vhat_H$ of the process $\vhat_{mH}$, that is, an 
orthogonal 
projection of the process $v_{H}$ parallel to the normalized score function $\beta_H$. Note that $\beta_H$ by its definition is of unit norm in the space 
$\Lcal_2(H)$. Moreover, functions $\beta_H$ and $\mathbb{1}$ are orthogonal. Therefore, from \eqref{em.pr-expand} we have the limit in 
distribution of the process $\vhat_{mH}(\phi)$ is
\begin{align}\label{limit-of-em.pr}
 \vhat_H (&\phi) = v_H(\phi) - \langle \phi, \beta_H \rangle_H v_H(\beta_H)\nonumber\\ 
 &\stackrel{\tiny{\mbox{from} (\ref{limitvH})}}{=} w_H(\phi) - \langle \mathbb{1}, \phi \rangle_H w_H(\mathbb{1}) - \langle \beta_H, \phi \rangle_H 
w_H(\beta_H).
\end{align}
This implies that $\vhat_H$ is an orthogonal projection of the function-parametric $H$-Brownian motion $w_H$ parallel to the subspace generated by 2 
functions $\set{\mathbb{1}, \beta_H}$. The process $\vhat_H$ depends not only on the true unknown exponent $\theta_0$ but also the score function 
$\beta_H$. In the terminology of \cite{Estate2}, the process $\vhat_H$ is called a $\beta_H$-projected H-Brownian motion.

\subsection{Method}\label{method}

Our approach follows the method in Khmaladze \cite{Estate2} for parametric family of distributions. The main idea can be briefly explained as follows: 
Under the null hypothesis, the empirical process $\vhat_{mH}(\phi)$ (see \eqref{empro} and \eqref{func.em.pr}) constructed on the tail of an unknown 
distribution $F$, from a fixed threshold $x_0$, possesses an ``unspecified'' limit in distribution (see \eqref{limit-of-em.pr}). ``Unspecified'' here 
means that the limit depends on some unknown parameter. We will map the 
process $\vhat_{mH}$ into another process $\vhat_{mG}$ (see \eqref{trans-em-pr-Gm}) whose limit in distribution is specified. The mapping procedure is 
one-to-one to guarantee that no statistical information is lost.

\medskip

We chose $G$ to be the exponential distribution $G(t)= 1-e^{-t}$ simply by preference with some reasons. That is, both $G$ and $H$ are members of 
the GPD family (see \eqref{gpd}). The distribution $G$ is the limiting distribution of the GPD$(k,\sigma)$ as $k \to 
0$ and scaled by $\sigma=1$. A research by Davison and Smith \cite{Davison} employed Kolmogorov-Smirnov and Anderson-Darling statistics to 
test the fit of the GPD to data using the critical value derived from an exponential distribution for the purpose. They discussed that the approach 
may be suspect since the exponential distribution is just a member of the whole family though.

\medskip

It is obvious that distributions $G$ 
and $H$ are equivalent or in other words, mutually absolutely continuous. It is also easy to check that the Fisher 
information of $G$ is $\Gamma_G =1$. Put
\[
\ell(t) = \sqrt{\frac{dG}{dH}(t)} = \theta^{-\frac{1}{2}} (1+t)^{\frac{\theta+1}{2}} e^{-\frac{t}{2}},
\]
then this function belongs to $\Lcal_2(H)$. In addition, if $\phi\in\Lcal_2(G)$ then $\ell\phi\in\Lcal_2(H)$ and $\norm{\phi}_G = \norm{\ell\phi}_H.$ Note 
that when 
$H$ and $G$ are equivalent, it is known that a transformation from a $H$-Brownian motion into a $G$-Brownian motion is straightforward by multiplication. 
Namely, $w_H(\ell \phi) = w_G(\phi)$ is a G-Brownian motion in $\Lcal_2(G)$. However, mapping a Brownian 
bridge such as $v_H$ or $\vhat_H$ to another Brownian bridge is not that straightforward any more. The fact is that $v_H(\ell\phi)$ depends on both $H$ 
and $G$ and so does $\vhat_H(\ell\phi)$.

\medskip

Consider a subspace $\Lcalhat$ of $\Lcal_2(H)$ generated by four functions $\set{\mathbb{1}, \beta_H, \ell,\ell\beta_G}$. These functions are of unit norm 
in 
$\Lcal_2(H)$. More explicitly, the score functions $\beta_H$ and $\beta_G$ are
\begin{align}\label{scorefuns}
\beta_H(t) = 1- \theta \log(1+t),\qquad \qquad \beta_G(t) = 1-t.
\end{align}

\noindent For any function $f$ and $g$ in $\Lcal_2(H)$, define the unitary operator $K_{f,g}$ as
\begin{align*}
K_{f,g} = I - \frac{1}{1-\langle f,g\rangle_H}(g-f) \langle g-f, \cdot \rangle_H
\end{align*}
where $I$ is the identity function. This unitary operator will turn the function $f$ into $g$ and reversely $g$ to $f$ and any function which is orthogonal 
to $f$ and $g$ in the subspace $\Lcal_2(H)$ into itself.

\medskip

Step into the method, we first consider the unitary operator
\begin{align}
K_{\mathbb{1},\ell} = I - \frac{1}{1-\langle \mathbb{1}, \ell\rangle_H} (\ell-\mathbb{1})\langle \ell-\mathbb{1}, \cdot \rangle_H.
\end{align}
This operator will map $\ell$ to $\mathbb{1}$ and $\mathbb{1}$ to $\ell.$ Next, consider the image of the function $\ell \beta_G$ via 
$K_{\mathbb{1},\ell}$, which is
\begin{align*}
\widetilde{\ell\beta}_G &= \ell \beta_G - \frac{1}{1-\langle \mathbb{1}, \ell\rangle_H}(\ell -\mathbb{1}) \langle \ell-\mathbb{1},\ell\beta_G\rangle_H\\
&= \ell \beta_G - \frac{1}{1-\int_0^{\infty} \ell(s)h(s)ds} (\ell -\mathbb{1})\int_0^{\infty} (\ell(s) -1) \ell(s) \beta_G(s) h(s) ds.
\end{align*}
Then consider the operator $K_{\beta_H,\widetilde{\ell\beta}_G}$ defined as 
\begin{align*}
K_{\beta_H,\widetilde{\ell\beta}_G} = I - \frac{1}{1-\langle \beta_H,\widetilde{\ell\beta}_G\rangle_H} (\widetilde{\ell\beta}_G -\beta_H) \langle 
\widetilde{\ell\beta}_G - \beta_H,\cdot\rangle_H.
\end{align*}

Set $\Khat = K_{\beta_H,\widetilde{\ell\beta}_G}K_{\mathbb{1},\ell}$, then this unitary operator will map $\ell$ to $\mathbb{1}$ and $\ell \beta_G$ to 
$\beta_H$. The non-uniqueness of such unitary operator like $\Khat$ was discussed thoroughly in Khmaladze \cite{Estate2}, Section 
3.4. Nevertheless, we believe that this operator $\Khat$ is simple enough for practical purpose, especially with only one parameter.

The main result for testing composite hypothesis was stated as Theorem 7 in \cite{Estate2}, it is essential so we restate it here accordingly to our 
notations.

\begin{theorem}\label{theorem} (Restatement of Theorem 7 in \cite{Estate2})

 If $\vhat_H$ is a $\beta_H$-projected $H$-Brownian motion and $G$ is absolutely continuous with respect to $H$, then
 \begin{align}
  \vhat_G(\phi) = \vhat_H(\Khat(\ell\phi)) = \Khat(\vhat_H(\ell\phi))
 \end{align}
is a $\beta_G$-projected G-Brownian motion.
\end{theorem}

\noindent As a consequence, transform the function-parametric tail empirical process $\vhat_{mH}(\phi)$ by $\Khat$, we obtain another process
\begin{align}\label{trans-em-pr-Gm}
 \vhat_{mG}(\phi) = \vhat_{mH}(\Khat(\ell\phi)) = \Khat(\vhat_{mH}(\ell\phi)),
\end{align}
which has $\vhat_G(\phi)$ as a limit in distribution.

\noindent Let 
\begin{align*}
\phi_x(t) = 1_{\set{t\le x}} 
\end{align*}
be series of indicator functions defined on $\Lcal_2(H)$ depending on $x$. Let $x$ runs from $0$ to $\infty$. From Theorem \ref{theorem}, we know that 
the 
limit in distribution $\vhat_H(\Khat(\ell\phi_x))$ of the process $\vhat_{mH}(\Khat(\ell\phi_x))$ is a $\beta_G$-projected G-Brownian motion in 
$\phi_x$. Hence, any statistic as an appropriate function based on $\vhat_{mH}(\Khat(\ell\phi_x))$ will be asymptotically distribution free. Denote by 
$\phitil_x$ the image of $\ell\phi_x$ via the operator $\Khat$, that is,
\begin{align*}
 \phitil_x = \Khat(\ell \phi_x) = K_{\beta_H,\widetilde{\ell\beta}_G}K_{\mathbb{1},\ell} (\ell\phi_x).
\end{align*}
For practical purpose, the form of $\phitil_x$ is
\begin{align}
 \phitil_x &= \ell \phi_x - \frac{1}{1-\langle \mathbb{1},\ell\rangle_H} (\ell-1) \langle \ell-1,\ell\phi_x\rangle_H - \frac{1}{1-\langle 
\beta_H,\widetilde{\ell\beta_G}\rangle_H}(\widetilde{\ell\beta_G}-\beta_H)\\\nonumber
 &\times \Big[\langle 
\widetilde{\ell\beta_G} -\beta_H,\ell\phi_x\rangle_H - \frac{\langle \ell-\mathbb{1},\ell\phi_x\rangle_H}{1-\langle 
\mathbb{1},\ell\rangle_H}\langle\widetilde{\ell\beta_G}-\beta_H,\ell-\mathbb{1}\rangle_H\Big].
\end{align}

\noindent Applying the unitary operator $\Khat$ on the process $\vhat_{mH}(\ell\phi_x)$ we have
\begin{align}\label{trans-em-pr-G}
 \vhat_{mG}(\phi_x) = \Khat(\vhat_{mH}(\ell\phi_x)) &= \vhat_{mH}(\phitil_x) = \int_0^{\infty} \phitil_x(t) \vhat_{mH}(dt)\nonumber\\
 &= \frac{1}{\sqrt{m}} \sum_{i=1}^m \big[\phitil_x(\Ttil_i) - \Esf_{\thetahat_m} \phitil_x(\Ttil_i)\big].
\end{align}
Here $\Esf_{\thetahat_m}$ denotes the expected value with respect to the distribution $H_{\thetahat_m}$. That means,
$\Esf_{\thetahat_m} f(\Ttil) = \int_0^{\infty} f(s) h_{\thetahat_m}(s) ds$ for any integrable function $f$. 
As a result of the Theorem \ref{theorem}, the limit in distribution of process $ \vhat_{mH}(\phitil_x)=\Khat(\vhat_{mH}(\ell\phi_x))$ is $\vhat_G(\phi_x)$ 
- a projected $G$-Brownian motion. We demonstrate in the next Section that any functionals of $\vhat_{mH}(\phitil_x)$ is asymptotically distribution 
free.

\section{Simulation results}\label{simulation}

The main purpose of this section is to show the asymptotically distribution free property of the new test statistics based on the transformed empirical 
process $\vhat_{mH}(\phitil_x)$. To do 
so, we will take some prevalent gof tests as examples, namely the Kolmogorov-Smirnov (KS), the Cramer-von Mises and the Anderson-Darling tests. These 
test statistics are some specific functionals of the transformed process $\vhat_{mH}(\phitil_x)$. Namely, an analogue version of the 
Kolmogorov-Smirnov test statistic is
\begin{align}
 KS = \max_{x}  \abs{\vhat_{mH}(\phitil_x)}.
\end{align}
The Cramer-von Mises statistics should be of the form
\begin{align}
 \Omega^2 = \int_0^{\infty} \vhat^2_{mH}(\phitil_x) dG(x),
\end{align}
and its weighted version called Anderson-Darling statistic is
\begin{align}
 A^2 = \int_0^{\infty} \frac{\vhat^2_{mH}(\phitil_x)}{G(x)(1-G(x))} dG(x).
\end{align}

\begin{figure}
	 \begin{subfigure}{.3\textwidth}
		\centering
		\includegraphics[width=1\linewidth]{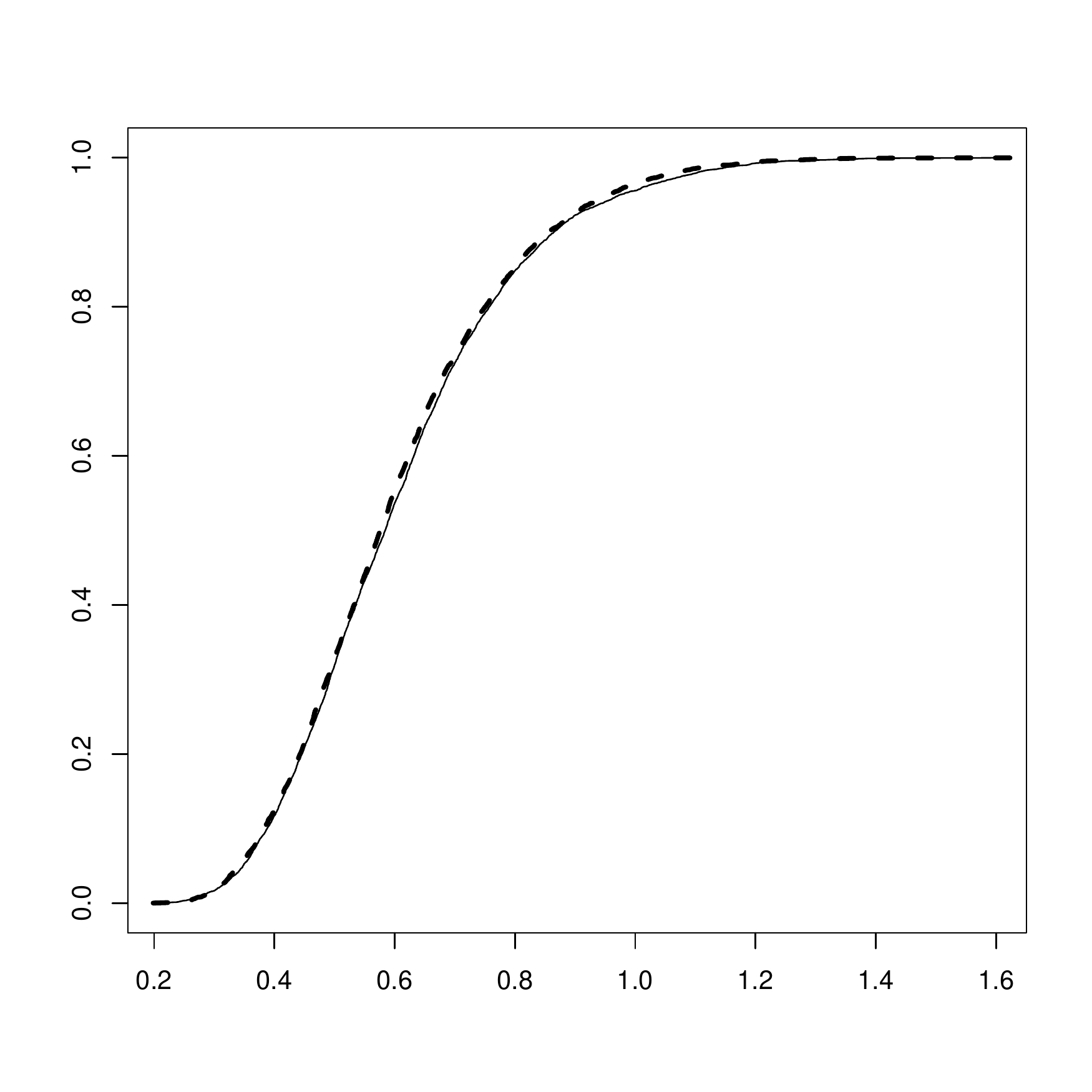}
		\subcaption{\it Threshold $x_0=3,$ sample size $n=1000$}
		\label{fig:x0-3-ss-1000}
	\end{subfigure}%
	\hspace{.03\linewidth}
	\begin{subfigure}{.3\textwidth}
		\centering
		\includegraphics[width=1\linewidth]{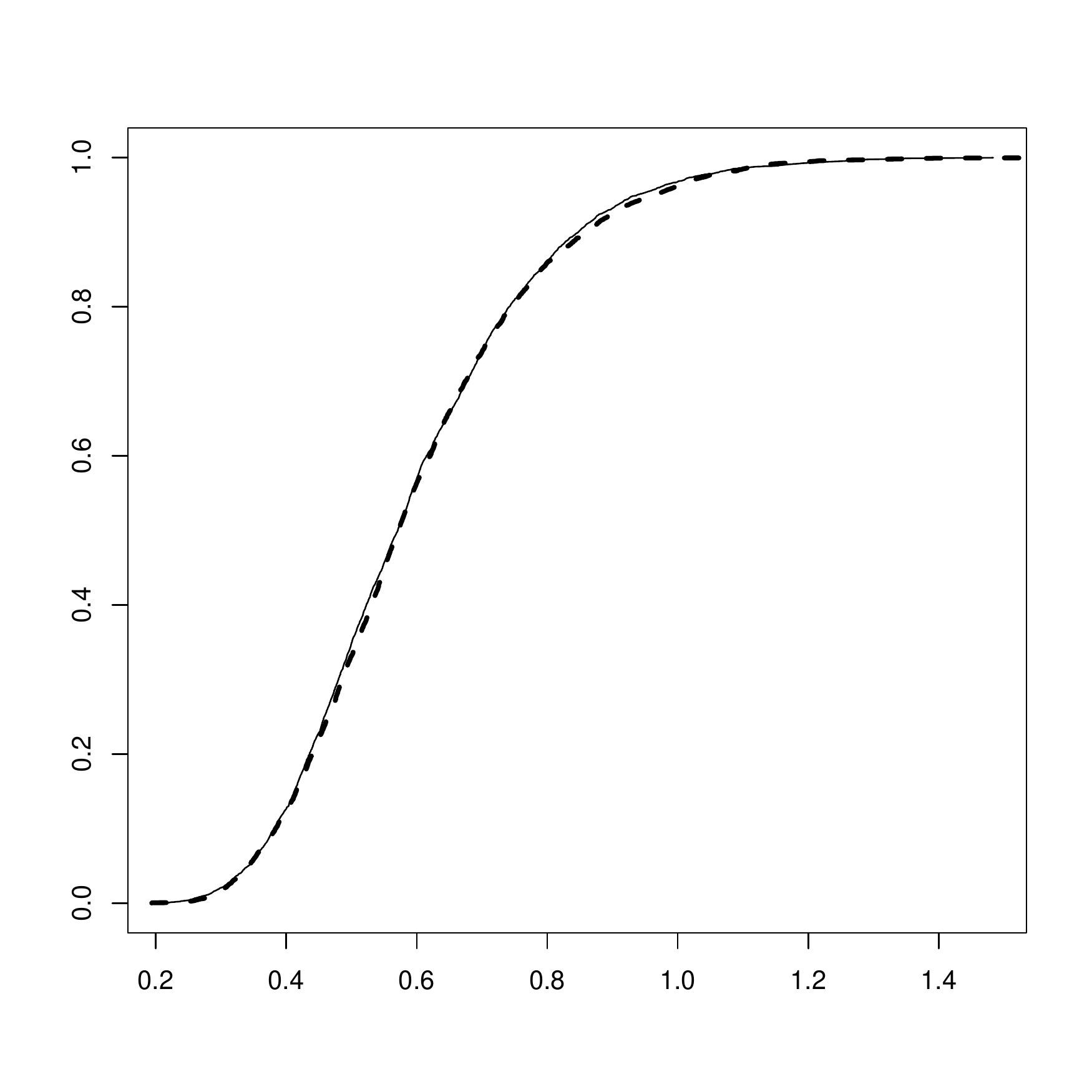}
		\subcaption{\it Threshold $x_0=5,$ sample size $n=5000$}
		\label{fig:x0-5-ss-5000}
	\end{subfigure}%
	\hspace{.03\linewidth}
	\begin{subfigure}{.3\textwidth}
		\centering
		\includegraphics[width=1\linewidth]{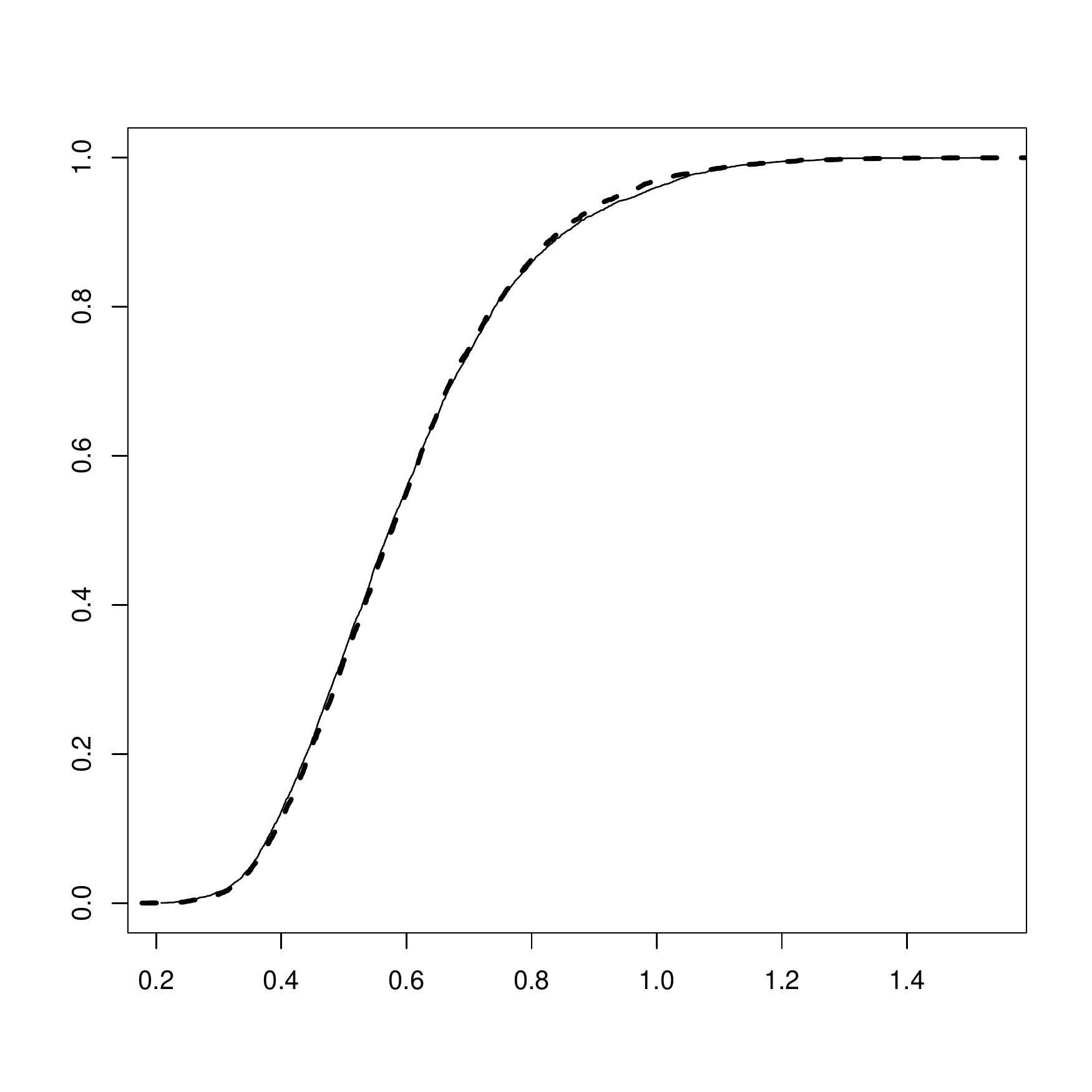}
		\subcaption{\it Threshold $x_0=10$, sample size $n=8000$}
		\label{fig:x0-10-ss-8000}
	\end{subfigure}
	\caption{\it Distribution of the KS test statistics. Solid line: Pareto distribution with $\theta_0=3$; Dashed line: Cauchy distribution}
	\label{fig:KS-test}
\end{figure}
\begin{figure}
	 \begin{subfigure}{.3\textwidth}
		\centering
		\includegraphics[width=1\linewidth]{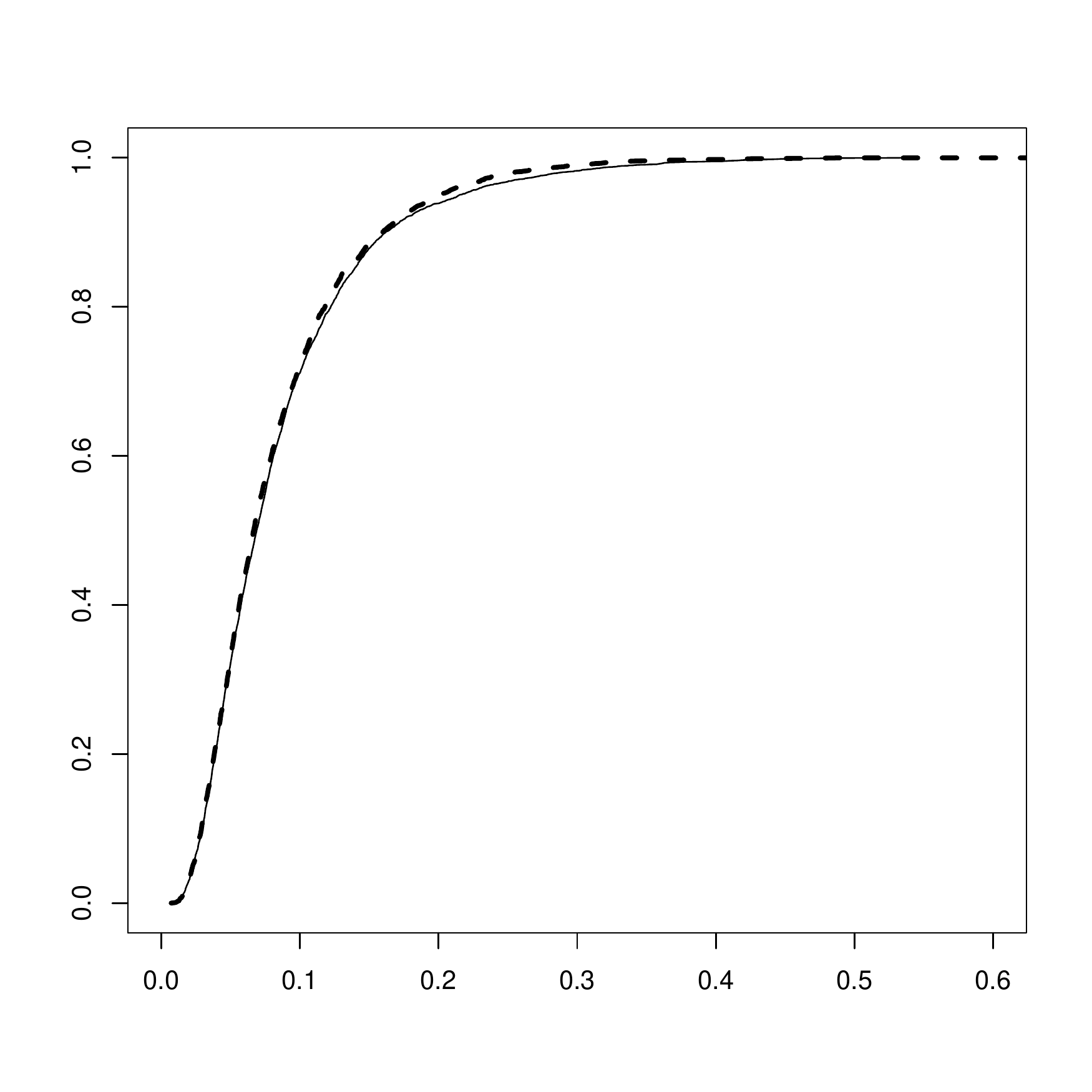}
		\subcaption{\it Threshold $x_0=3,$ sample size $n=1000$}
		\label{fig:omega-x0-3-ss-1000}
	\end{subfigure}%
	\hspace{.03\linewidth}
	\begin{subfigure}{.3\textwidth}
		\centering
		\includegraphics[width=1\linewidth]{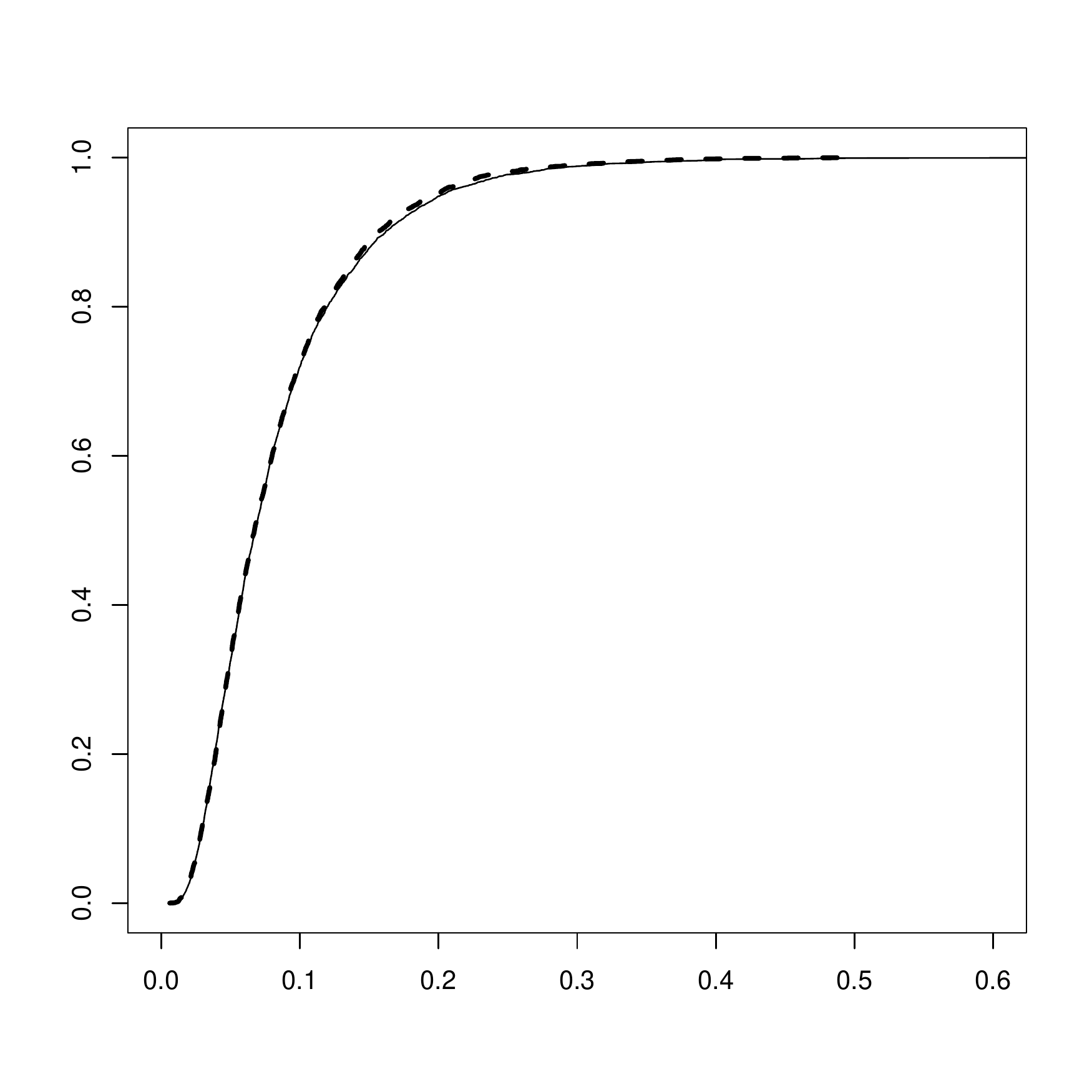}
		\subcaption{\it Threshold $x_0=5,$ sample size $n=5000$}
		\label{fig:omega-x0-5-ss-5000}
	\end{subfigure}%
	\hspace{.03\linewidth}
	\begin{subfigure}{.3\textwidth}
		\centering
		\includegraphics[width=1\linewidth]{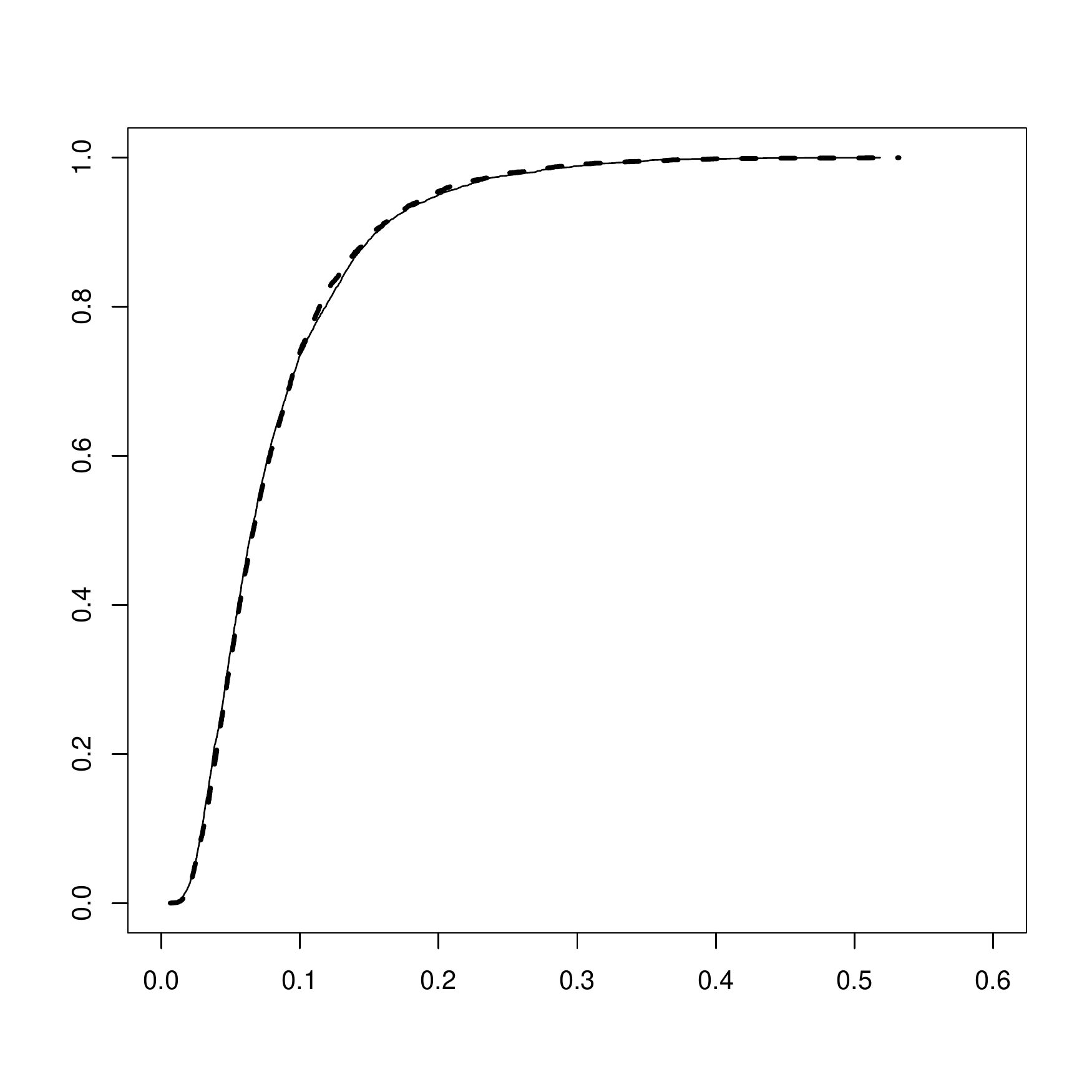}
		\subcaption{\it Threshold $x_0=10$, sample size $n=5000$}
		\label{fig:omega-x0-10-ss-5000}
	\end{subfigure}
	\caption{\it Distribution of the $\Omega^2$ test statistics. Solid line: Pareto distribution with $\theta_0=2$; Dashed line: Cauchy distribution}
	\label{fig:omega-test}
\end{figure}
\begin{figure}
	 \begin{subfigure}{.3\textwidth}
		\centering
		\includegraphics[width=1\linewidth]{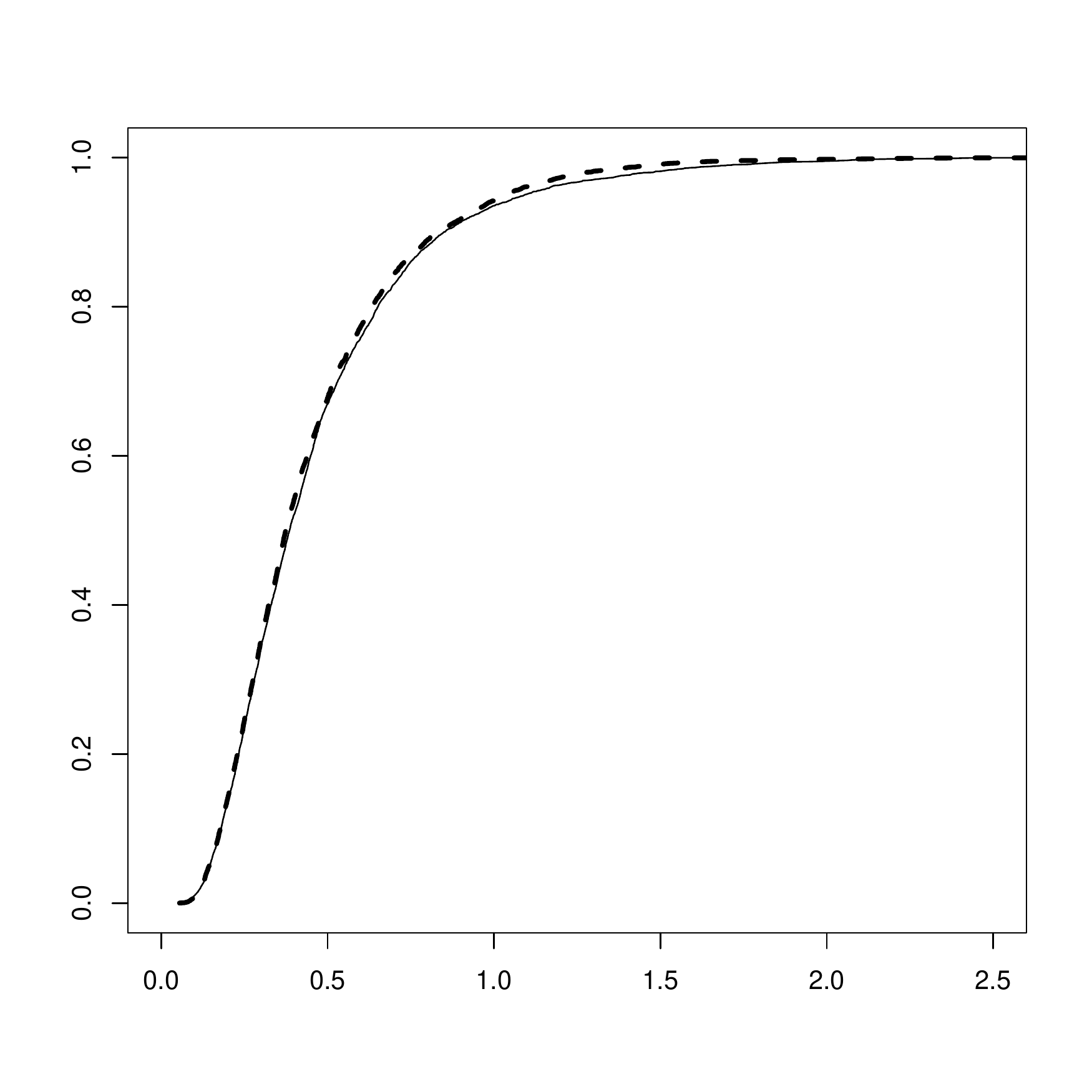}
		\subcaption{\it Threshold $x_0=3,$ sample size $n=1000$}
		\label{fig:omega-weighted-x0-3-ss-1000}
	\end{subfigure}%
	\hspace{.03\linewidth}
	\begin{subfigure}{.3\textwidth}
		\centering
		\includegraphics[width=1\linewidth]{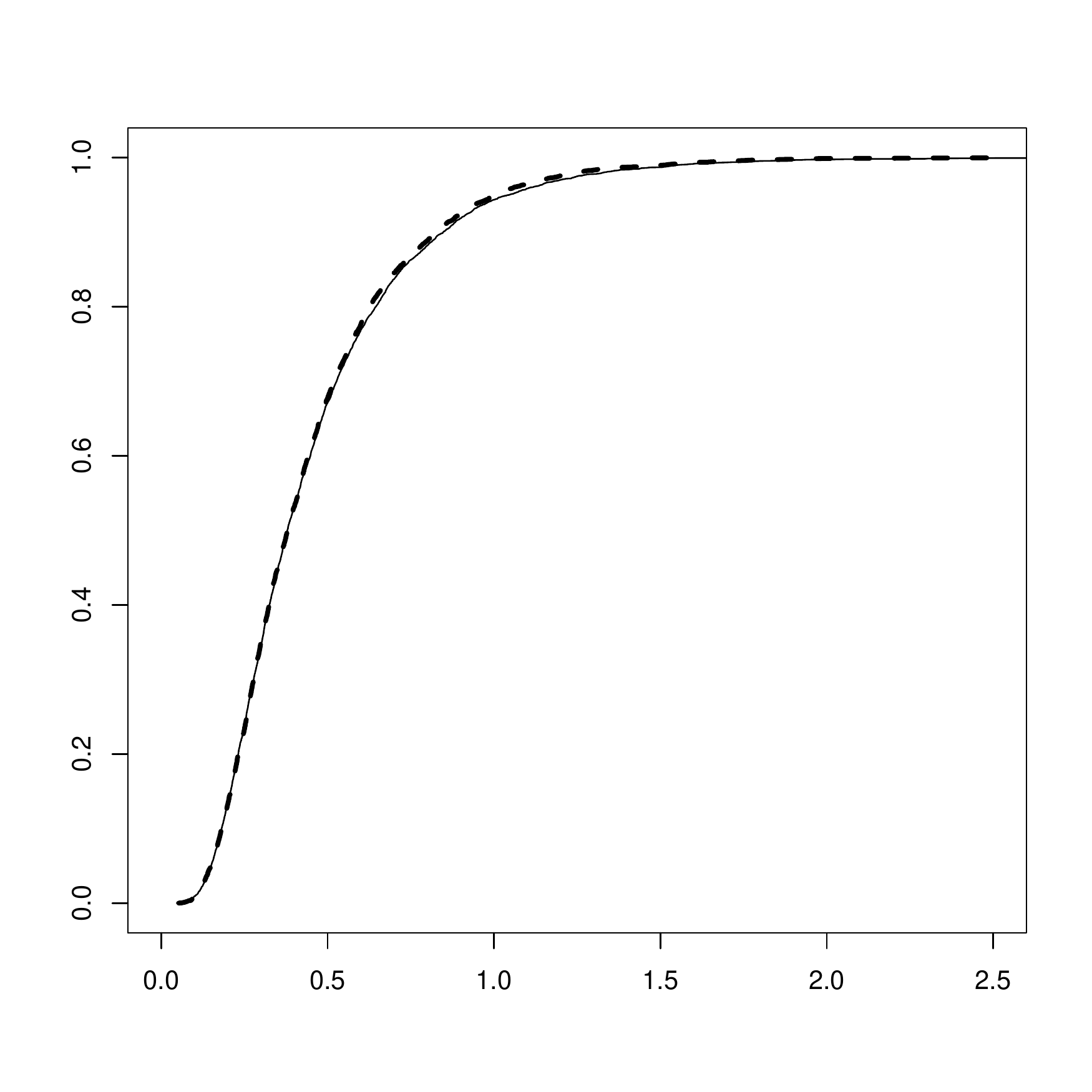}
		\subcaption{\it Threshold $x_0=5,$ sample size $n=5000$}
		\label{fig:omega-weighted-x0-5-ss-5000}
	\end{subfigure}%
	\hspace{.03\linewidth}
	\begin{subfigure}{.3\textwidth}
		\centering
		\includegraphics[width=1\linewidth]{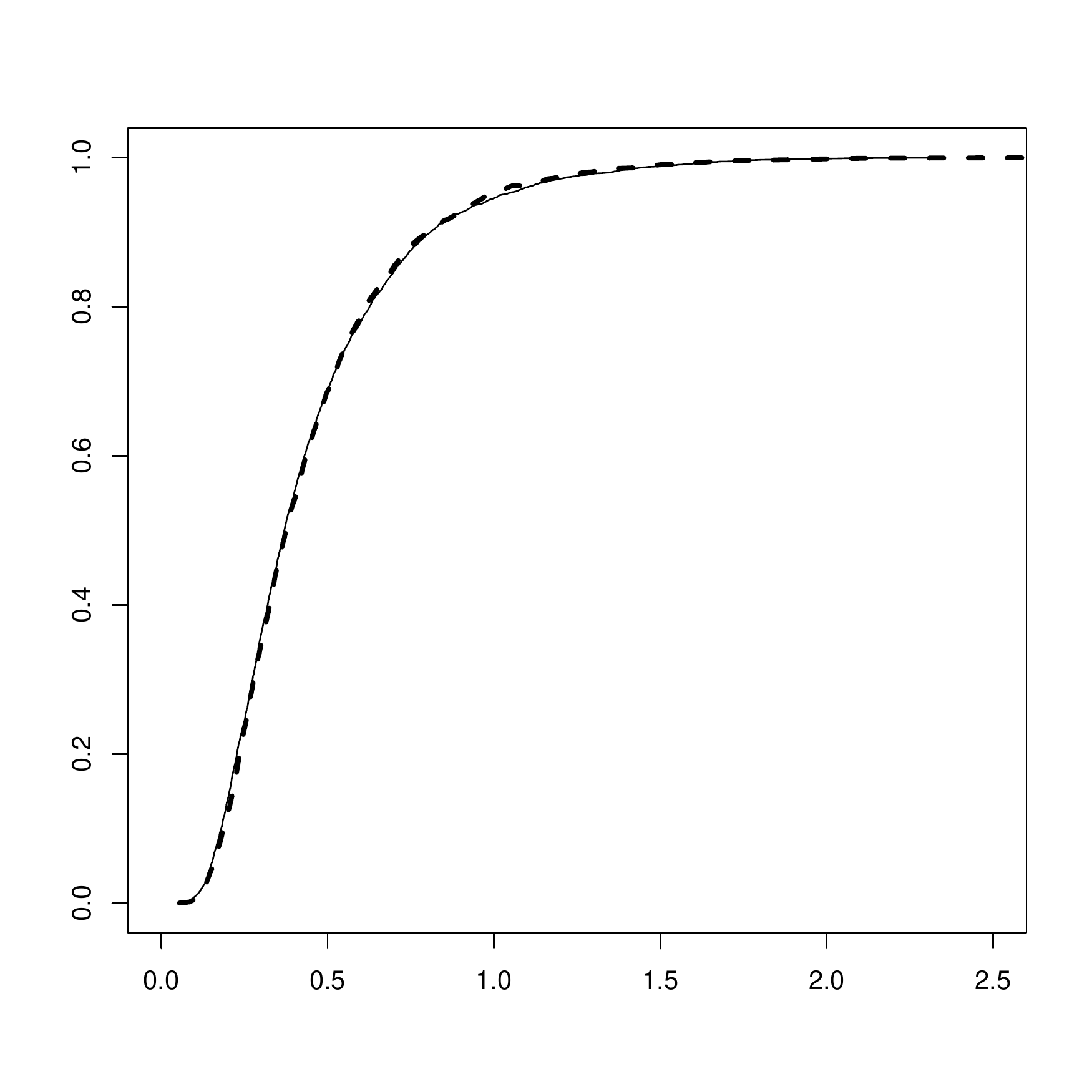}
		\subcaption{\it Threshold $x_0=10$, sample size $n=5000$}
		\label{fig:omega-weighted-x0-10-ss-5000}
	\end{subfigure}
	\caption{\it Distribution of the $A^2$ test statistics. Solid line: Pareto distribution with $\theta_0=0.5$; Dashed line: Cauchy 
distribution}
	\label{fig:anderson-darling-test}
\end{figure}

\newpage
As shown in Section \ref{method}, the limit in distribution of the process $\vhat_{mH}(\phitil_x)$ is a projected Brownian motion in time 
$G(x)$. Hence, theoretically these test statistics are asymptotically distribution free. Especially, it makes sense that we integrate with respect to 
nothing else but $G(x)$ in order to get the $\Omega^2$ and $A^2$ tests. In principle, $x$ runs from $0$ to infinity. However, we choose $x\in 
\set{0.1,0.2, \cdots, 7.9,8}$ since $G(8) = 0.9997 \approx 1$ hence we practically just need $x_{\max}=8$ as the maximum value for $x$. The $\Omega^2$ and 
$A^2$ tests can be approximated as follows:
\begin{align}
 \Omega^2 &\approx 0.1 \times \sum_{x=0.1}^{8}  \vhat_{mH}(\phitil_x) e^{-x},\nonumber\\
 A^2 &\approx 0.1 \times  \sum_{x=0.1}^{8}  \frac{\vhat_{mH}(\phitil_x)}{1-e^{-x}}.\nonumber
\end{align}

To create the curves of the cumulative distributions of the new tests, we choose Pareto and Cauchy distributions as the underlying distributions. It is 
known that the tail of the Cauchy distribution is regularly varying with the exponent $\theta_0 =1$. For the Pareto distributions, the exponent is 
arbitrarily selected and positive. We did choose some $\theta_0$ ranging from $0.5$ to $10$. Sample size $n$ for simulation needs to be large to guarantee 
that the tail is sufficiently big, namely, $m$ is not less than $40$. Hence, we chose $n$ at least equals $1000$. We also chose different threshold $x_0$, 
namely $x_0 =3,5,10$. Usually, each curve is produced by $5000$ iterations of simulation. 

Figures \ref{fig:KS-test}, \ref{fig:omega-test} and \ref{fig:anderson-darling-test} show the plots of the cumulative functions for $KS, \Omega^2, A^2$ 
test statistics respectively with different choices of the threshold $x_0$, different sample sizes and different underlying distributions. As we can see 
in these figures, the two curves of the cumulative 
distribution functions in each plot are not distinguishable, which practically demonstrate the asymptotically distribution free property of the new test 
statistics. 
Moreover, as we notice, for different $x_0$, the difference between these curves is also very minor.

Regarding the time of the procedure, it took approximately $1$ hour to create the cumulative distribution functions with $x_0 =5$ and sample size $n=5000$ 
by $5000$ iterations for two different original distributions $F$. Therefore, we believe that the method is easy to implement and also quick.

\section{Discussion}

It is possible to extend this approach to testing multidimensional regularly varying tail distributions. As long as we are able to find a transparent expression for the class of regularly varying tail distributions with an explicit parametric family form, we are able to transform the whole family into a specific distribution, from that we can build up a whole class of asymptotic distribution-free test statistics. 

\bigskip

\noindent {\bf Acknowledgements}

The author is greatly indebted to Professor Estate Khmaladze for his invaluable suggestion for approaching this 
problem 
and many useful discussions during this research. This research was conducted during the author's doctorate of philosophy degree, so the author would also like to thank Victoria University of Wellington for providing the doctoral scholarship to conduct this research.





\end{document}